\documentclass[12pt,twoside]{amsart}
\usepackage{amsopn,amssymb}
\usepackage[margin=1in]{geometry}
\nonstopmode

\DeclareMathOperator{\link}{link}
\DeclareMathOperator{\fdel}{del}

\newcommand{\De}{\Delta}

\newcommand{\si}{\sigma}
\newcommand{\ta}{\tau}
\newcommand{\QQ}{\mathbb{Q}}

\newcommand{\st}{\; | \;}

\theoremstyle{definition}
\newtheorem*{acknowledgement}{Acknowledgement}

\begin{document}

\title{S\MakeLowercase{implicial}D\MakeLowercase{ecomposability}: \MakeLowercase{a package for }M\MakeLowercase{acaulay} 2}

\author[D.\ Cook II]{David Cook II}
\address{Department of Mathematics, University of Kentucky, 715 Patterson Office Tower, Lexington, KY 40506-0027, USA}
\email{dcook@ms.uky.edu}
\thanks{This article describes version 1.0.4 of SimplicialDecomposability.  The current version can be downloaded at {\em http://www.ms.uky.edu/$\sim$dcook/files/SimplicialDecomposability.m2}.}
\date{\today}

\begin{abstract}
    We introduce a new {\em Macaulay 2} package, {\em SimplicialDecomposability}, which works in conjunction with the extant package
    {\em SimplicialComplexes} in order to compute a shelling order, if one exists, of a specified simplicial complex.  Further,
    methods for determining vertex-decomposability are implemented, along with methods for determining $k$-decomposability.
\end{abstract}

\maketitle

\subsection*{Introduction}

A {\em simplicial complex} $\De$ on a finite vertex set $V$ is a set of subsets of $V$ closed under inclusion.  Elements $\si \in \De$
are called {\em faces} and the {\em dimension} of $\si$ is $\#\si - 1$.  The {\em dimension} of $\De$ is $\max\dim\si$.  The
{\em $f$-vector} of a $\De$, where $d = \dim\De + 1$, is the $(d+1)$-tuple  $(f_{-1}, \ldots, f_{d-1})$, where $f_i$ is the number
of faces of dimension $i$ in $\De$.  Using this, the {\em $h$-vector} of $\De$ is the $d+1$-tuple $(h_0, \ldots, h_d)$ given by
$h_j = \sum_{i=0}^j(-1)^{j-i}\binom{d-i}{j-i}f_{i-1}$ for $0 \leq j \leq d$.  

The {\em Stanley-Reisner ideal} is the ideal $I(\De)$ generated by the minimal non-faces of $\De$ and the {\em Stanley-Reisner ring}
is the ring $K[\De] = K[V]/I(\De)$, for a given field $K$.  Thus the Stanley-Reisner ideals of complexes on a given vertex set $V$ are
exactly the squarefree monomial ideals in $K[V]$.  Using relations between the complex and the ideal, one can use tools from both
algebra and combinatorics to study properties of both.  For example, the $h$-vector of a complex $\De$ is the coefficient-vector of
the numerator of the Hilbert series of $K[\De]$.

The package {\em SimplicialComplexes} by Sorin Popescu, Gregory G. Smith, and Mike Stillman already implements many methods for 
simplicial complexes in {\em Macaulay 2}~\cite{M2}, a software system designed to aid in research of commutative algebra and algebraic
geometry.  We introduce a new package, {\em SimplicialDecomposability}, for {\em Macaulay 2} which provides several new methods for
testing various forms of decomposability for simplicial complexes.  Particularly, the package implements methods for testing shellability
and vertex-decomposability.

\subsection*{Shellability}

If a simplicial complex has one facet, say $\si$, then it is a simplex and is denoted $2^\si$.  Let $\De$ be a simplicial complex which
has equi-dimensional facets, i.e., is {\em pure}.  Then by \cite[Definition~III.2.1]{St}, $\De$ is {\em shellable} if its facets can be
ordered $\si_1, \ldots, \si_n$ such that
\[
    \bigcup_{j=1}^{i}{2^{\si_j}} \setminus \bigcup_{j=1}^{i-1}{2^{\si_j}}
\]
has a unique minimal element for $2 \leq i \leq n$, such an ordering is called a {\em shelling order}.

See \cite[Definition~2.1]{BW-1} for the definition of non-pure shellability, which is implemented in the package for non-pure complexes.

Shellability is of interest because it implies a number of nice properties.  In particular, if a pure simplicial complex is shellable,
then its Stanley-Reisner ring is Cohen-Macaulay over every field~\cite[Theorem~III.2.5]{St}.  Hence, its $h$-vector is
non-negative and can be read off from any shelling order~\cite[Theorem~III.2.3]{St}.  Further still, the $h$-vectors of pure
shellable complexes are numerically classified~\cite[Theorems~II.2.2 and~II.3.3]{St}.

We recall that the {\em Alexander dual} of a simplicial complex $\Delta$ on vertex set $V$ is the simplicial complex
$\De^\vee := \{V \setminus F \st F \not\in \Delta\}$.  Further, we say an ideal $I = (f_1, \ldots, f_n)$ has {\em linear quotients}
if for $1 < i \leq n$, the quotient ideal $(f_1, \ldots, f_{i-1}):(f_i)$ is generated by linear forms.

In the following example we demonstrate~\cite[Theorem~1.4(c)]{HHZ} which shows that a pure simplicial complex is shellable if and only if
the Stanley Reisner ideal of the Alexander dual has linear quotients.  We begin by constructing the polynomial ring $R = \QQ[a,b,c,d,e,f,g]$
and a simplicial complex $D$, which we verify is pure.  Note that loading the package {\tt SimplicialDecomposability} automatically loads
the package {\tt SimplicialComplexes}.

\begin{verbatim}
  i1 : needsPackage "SimplicialDecomposability";
  i2 : R = QQ[a..g];
  i3 : D = simplicialComplex monomialIdeal {a*b,a*c,b*c,c*d,d*e,d*f,f*g};
  i4 : isPure D
  o4 = true
\end{verbatim}

We can recover the sequence of linear quotients directly from a shelling order.  We recall that a pure simplicial
complex $\De$ is shellable if there is an order of the facets $F_1, \ldots, F_n$ such that for $0 < j < i$ there exists an
$x \in F_i \setminus F_j$ and a $0 < k < i$ such that $F_i \setminus F_k = \{x\}$.  The set of vertices associated to each $i$
in the preceding statement generate the linear quotient order of $I(\De^\vee)$ with respect to the given shelling order (see
the proof of~\cite[Theorem~1.4(c)]{HHZ}).

\begin{verbatim}
  i5 : -- find the linear quotients from a shelling order
       linearQuotients = O -> (
         for i from 1 to #O-1 list
           unique flatten for j from 0 to i-1 list (
             ImJ = set support O_i - set support O_j;
             for k from 0 to i - 1 list (
               ImK = set support O_i - set support O_k;
               if #ImK == 1 and isSubset(ImK, ImJ) then
                 first toList ImK else continue)));
\end{verbatim}

We generate a shelling order $O_1$ of $D$ with the method {\tt shellingOrder}.  This method attempts to build up a shelling order of $D$
recursively using a depth-first search, adding one facet at a time.  We note that in the non-pure case, the method only searches the
remaining facets of largest dimension.

\begin{verbatim}
  i6 : O1 = shellingOrder D
  o6 = {c*e*g, b*e*g, a*e*g, b*d*g, a*d*g, c*e*f, b*e*f, a*e*f}
  o6 : List
  i7 : linearQuotients O1
  o7 = {{b}, {a}, {d}, {d, a}, {f}, {f, b}, {f, a}}
  o7 : List
\end{verbatim}

It is sometimes beneficial to have more than one shelling order for a given simplicial complex.  We can use the option {\tt Random}
with the method {\tt shellingOrder} to first apply a random permutation to the facets before preceding with the recursion.

\begin{verbatim}
  i8 : O2 = shellingOrder(D, Random => true)
  o8 = {b*e*g, a*e*g, a*d*g, b*e*f, c*e*f, a*e*f, c*e*g, b*d*g}
  o8 : List
  i9 : linearQuotients O2
  o9 = {{a}, {d}, {f}, {c}, {f, a}, {c, g}, {d, b}}
  o9 : List
\end{verbatim}

Alternately, we may use the option {\tt Permutation} with the method {\tt shellingOrder} to force a given permutation on the facets
before preceding with the recursion.

\begin{verbatim}
  i10 : O3 = shellingOrder(D, Permutation => {3,2,1,0,4,5,6,7})
  o10 = {b*d*g, b*e*g, a*e*g, c*e*g, a*d*g, c*e*f, b*e*f, a*e*f}
  o10 : List
  i11 : linearQuotients O3
  o11 = {{e}, {a}, {c}, {a, d}, {f}, {f, b}, {f, a}}
  o11 : List
\end{verbatim}

Thus we now have multiple distinct linear quotient orders associated to the ideal $I(D^\vee)$, each coming from a distinct
shelling order of $D$.

\subsection*{Vertex-decomposability}

Let $\De$ be a pure simplicial complex and $\si$ a face of $\De$.  Then the {\em link} and {\em face deletion} of $\De$ by $\si$ are
the simplicial complexes
\[
    \link_{\De}{\si} := \{ \ta \in \De \st \si \cap \ta = \emptyset, \si \cup \ta \in \De \}
    \mbox{ and }
    \fdel_{\De}{\si} := \{ \ta \in \De \st \si \nsubseteq \ta \}.
\]
Then \cite[Definition~2.1]{PB} defines $\De$ to be {\em vertex-decomposable} if either $\De$ is a simplex or there exists a vertex
$x \in \De$, called a {\em shedding vertex}, such that $\link_{\De}{x}$ and $\fdel_{\De}{x}$ are vertex-decomposable.  

See \cite[Definition~11.1]{BW-2} for the definition of non-pure vertex-decomposability, which is implemented in the package for non-pure complexes.
Also, see~\cite[Definitions~3.1 and~3.6]{Wo} for the generalisation of vertex-decomposability, called $k$-decomposability.  It is implemented in
the package with the methods {\tt iskDecomposable} and {\tt isSheddingFace}.

Being vertex-decomposable is a strong property which implies many things.  A pure vertex-decomposable simplicial complex is
shellable~\cite[Theorem~2.8]{PB} and hence has non-negative $h$-vector~\cite[Theorem~III.2.3]{St} and its Stanley-Reisner ring is
Cohen-Macaulay~\cite[Theorem~III.2.5]{St}.  Furthermore, the $h$-vectors are numerically classified for vertex-decomposable
simplicial complexes~\cite[Theorem~3.5]{Lee}.  Moreover, the Stanley-Reisner ring of a pure vertex-decomposable complex is
squarefree glicci, that is, in the {\bf G}orenstein {\bf li}aison {\bf c}lass of a {\bf c}omplete {\bf i}ntersection such that the
even links are squarefree monomials \cite[Theorem~3.3]{NR}.

In the following example we demonstrate that the simplicial complex $D$ from the previous example is indeed squarefree glicci.  We
use~\cite[Remark~2.4]{NR} to find a basic double link of $I(D)$ to $I(\link_D{v})$, both in $R$, for some shedding vertex $v$ of $D$.  

First, we verify that $D$ is vertex-decomposable.  We then find all of its shedding vertices.

\begin{verbatim}
  i12 : isVertexDecomposable D
  o12 = true
  i13 : select(allFaces(D, 0), v -> isSheddingVertex(D, v))
  o13 = {a, b, c, d, f}
  o13 : List
\end{verbatim}

We choose the shedding vertex $f$ of $D$ and generate $E = \link_D{f}$.  Moreover, we in turn find its shedding vertices.

\begin{verbatim}
  i14 : E = link(D, f); ideal E
  o15 = ideal (a*b, a*c, b*c, d, f, g)
  o15 : Ideal of R
  i16 : select(allFaces(E, 0), v -> isSheddingVertex(E, v))
  o16 = {a, b, c}
  o16 : List
\end{verbatim}

We now choose the shedding vertex $c$ of $D$ and generate $F = \link_E{c}$.  Notice then that the Stanley Reisner ideal of 
$F$ is a complete intersection.

\begin{verbatim}
  i17 : F = link(E, c); ideal F
  o18 = ideal (a, b, c, d, f, g)
  o18 : Ideal of R
\end{verbatim}

Hence, we now have the following sequence of basic double links in $R$ which has squarefree terms on the even steps (the
odd steps are omitted):
\[
    \QQ[D] = (ab, ac, bc, cd, de, df, fg) \sim \QQ[E] = (ab, ac, bc, d, f, g) \sim \QQ[F] = (a,b,c,d,f,g).
\]

%

\begin{acknowledgement}
    The author would like to thank his advisor, Uwe Nagel, for reading drafts of this article and making comments thereover.  The author would
    also like to thank Russ Woodroofe for pointing out the non-pure generalisation of $k$-decomposability in~\cite{Wo} and an anonymous
    referee for inspiring a more concrete example herein.

    Part of the work for this paper was done while the author was at a Macaulay 2 workshop in Berkeley, California, January 8, 2010 through
    January 12, 2010, organized by Amelia Taylor and Hirotachi Abo with David Eisenbud, Daniel R. Grayson, and Michael E. Stillman, and
    funded by the National Security Agency (NSA) through grant H98230-09-1-0111.
\end{acknowledgement}


\end{document}